\documentclass [a4paper, 12pt, fleqn] {article}

\usepackage{amssymb}
\usepackage{amsmath}
\usepackage{amsthm}
\usepackage[all]{xy}
\newcommand{\Lb}[1]{\mathcal{L}_{#1}}
\newcommand{\fa}{\ \forall \,}
\newcommand{\dpp}{\mathrel{\mathop:}=}
\newcommand{\N}{\mathbb{N}}
\renewcommand{\P}{\mathbb{P}}
\newcommand{\oo}{\mathcal{O}}
\newcommand{\M}{\mathcal{M}}
\newcommand{\I}{\mathcal{I}}
\renewcommand{\S}{\mathcal{S}}
\newcommand{\ro}{\varrho}
\newcommand{\End}{{\rm{End}}}
\newcommand{\Hom}{{\rm{Hom}}}

\renewcommand{\Im}{{\rm{Im}}}
\newcommand{\Pic}{{\rm{Pic}}}
\newcommand{\gr}{{\rm{gr}}}
\newcommand{\res}{{\rm{res}}}
\newcommand{\car}{{\rm{char}}}
\newcommand{\id}{{\rm{id}}}
\newcommand{\supp}{{\rm{supp}}}

\newcommand{\lex}{{\rm{lex}}}

\newcommand{\ad}{{\rm{ad}}}
\newcommand{\diag}{{\rm{diag}}}
\newcommand{\ueber}{\!\!\mid\!}
\newcommand{\betr}[1]{|#1|}
\newcommand{\rb}[1][1.5]{\raisebox{-#1ex}{}}

\newtheoremstyle{normal}{2,5ex minus1ex}{2,5ex minus1ex}{\upshape}{}{\bfseries}{.}{ }{}
\newtheoremstyle{bemerkung}{\parskip}{\parskip}{\upshape}{}{\bfseries}{.}{}{}
\theoremstyle{normal}
\newtheorem{satz}{Proposition}
\newtheorem{lemma}[satz]{Lemma}
\newtheorem{folg}[satz]{Corollary}
\newtheorem{prop}[satz]{Proposition}
\newtheorem{defi}{Definition}
\newtheorem*{notation}{Notation}
\newtheorem{bem}{Remark}
\newenvironment{beweis}{\begin{proof}[\textbf{\upshape Proof.\\}]}{\end{proof}}
\newcommand{\exendbeweis}{\renewcommand{\qed}{\relax}\end{beweis}}
\newcounter{aufzcounter}
\newenvironment{aufz}{
\begin{list}{\arabic{aufzcounter})}{\usecounter{aufzcounter} \itemsep0.2ex \topsep0.2ex \parsep\parskip \labelwidth2ex \leftmargin5ex}}
{\end{list}}

\begin{document}

\title{Standard monomials for wonderful group compactifications} 
\author{Katrin Appel}
\date{}
\maketitle

\begin{abstract}
Let $X$ be the wonderful compactification of the semisimple adjoint algebraic group $G$. We show that the basis of $H^0(X,\Lb{})$ constructed in \cite{CM03} is compatible with all $B \times B$-orbit closures in $X$ by defining subsets using only combinatorics of the underlying paths. Furthermore, we construct standard monomials on $X$ that have properties similar to classical standard monomials. 
\end{abstract}


\section{Introduction}
 
Let $G$ be a semisimple algebraic group over an algebraically closed field. One aim of standard monomial theory is to define bases of the space of global sections $H^0(G/B,\Lb\lambda)$. Here $B$ is a Borel subgroup of $G$, $G/B$ is the flag variety, and $\lambda$ is a dominant weight of $G$ with associated line bundle $\Lb\lambda$. The elements of this basis should be weight vectors of $G$ and behave nicely under restriction to the Schubert varieties in $G/B$. As $H^0(G/B,\Lb\lambda) \cong V(\lambda)^*$ by the Borel-Weil theorem, this also gives bases of the highest weight modules of $G$ that have nice geometric properties. 

One solution for this is the path model which defines paths in the weight lattice and an associated path vector $p_\pi \in H^0(G/B,\Lb\lambda)$ for every LS-path of the form $\lambda$. These path vectors form a basis of the $G$-module $H^0(G/B,\Lb\lambda)$. Given a Schubert variety $S(w)$ in $G/B$, it is possible to define a subset of the set of all LS-paths of the form $\lambda$ such that the associated path vectors restrict to a basis of $H^0(S(w),\Lb\lambda)$. This definition depends purely on combinatorial properties of the paths. The path vectors associated to LS-paths are called standard monomials an $G/B$.

As a generalization, it would be nice to have an analogue for compactifications of symmetric spaces $G/H$ instead of the flag variety $G/B$. This paper deals with the special case of the wonderful compactification of an adjoint group $G$, that is considered as a symmetric space of the group $G \times G$. The wonderful compactification $X$ of the group $G$ consists of several $G \times G$-orbits of which exactly one is closed, and this unique closed orbit $Y$ is isomorphic to $G/B \times G/B$. Now the existing standard monomials on $G/B \times G/B$ can be extended to $X$ to get standard monomials on $X$.

In this paper we show, that arbitrary extensions are compatible not only with the closures of $G \times G$-orbits in $X$, but also with the closures of all $B \times B$-orbits. Furthermore, we construct extensions that possess attributes similar to those of classical standard monomials.

\vspace{1ex}

The paper is organized as follows.
In the second section, we shortly recall the construction of the wonderful compactification $X$ of a group $G$. We also include the description of the $B \times B$-orbits in $X$ and their closure relations obtained by Springer. Furthermore, a short description of the line bundles on $X$ is given.

In the third section, we very briefly describe standard monomials on the flag variety $G/B$. We do not give an adequate presentation of the path model and standard monomials, but state only those properties of LS-paths and path vectors we need later.

The fourth section contains a resume of the results of the paper \cite{CM03}, in which Chiriv\`{i} and Maffei construct standard monomials for wonderful compactifications of symmetric spaces. We give the definition and the properties in the special case of the wonderful compactification of a group. In particular, the set $\M^{(\lambda)}$ is defined which is a basis for $H^0(X,\Lb\lambda)$ compatible with restrictions to $G \times G$-orbit closures.

In the fifth section, we define for every $B \times B$-orbit closure $Z$ a subset $\M^{(\lambda)}_Z$ of $\M^{(\lambda)}$ which restricts to a linearly independent subset of $H^0(Z,\Lb\lambda)$. Using consequences of the existence of a compatible Frobenius splitting and properties of the associated graded module for a suitable filtration, we prove that $\M^{(\lambda)}_Z$ is also a basis of the global sections on $Z$. Thus, we show that in the case of a group compactification the basis $\M^{(\lambda)}$ is compatible with restrictions to $B \times B$-orbit closures as well.

In the last section we take advantage of the fact that so far, none of the constructions depends on the choice of the continuation of the standard monomials on $Y$ to $X$. This section is dedicated to the construction of standard monomials for $X$ which have properties similar to the classical standard monomials.

\vspace{1ex}

The author wants to thank P.~Littelmann and M.~Brion for useful comments.


\section{Wonderful group compactifications}

In \cite{CP83}, De Concini and Procesi construct the wonderful compactification of an adjoint group $G$ in characteristic zero as a special case of symmetric spaces. The definition was extended to positive characteristic by Strickland in \cite{St87}. In this section, we shortly recall the construction and properties of wonderful group compactifications.

Let $G$ be an adjoint semisimple algebraic group over an algebraically closed field of arbitrary characteristic. Choose a Borel subgroup $B$ and a torus $T$ such that $T \subset B \subset G$. The corresponding weight lattice is $\Lambda$, with dominant weights $\Lambda^+$. The set of simple roots is $\Delta=\{\alpha_1,\ldots,\alpha_l\}$, the positive roots are denoted by $\Phi^+$ and the negative roots by $\Phi^-=-\Phi^+$. 
Let $\tilde{G}$ be the simply connected covering of $G$ with morphism $\pi_{\ad}: \tilde{G} \rightarrow G$. Denote the pre-images of the subgroups of $G$ in $\tilde{G}$ by $\tilde{T}=\pi_{\ad}^{-1}(T)$, $\tilde{B}=\pi_{\ad}^{-1}(B)$. 

The group $G$ can be considered as the symmetric space $(G \times G)/\diag_G$, where $\diag_G$ denotes the diagonal in $G \times G$. For a suitable $G$-module $M$, this space $G \cong (G \times G)/\diag_G$ is isomorphic to the orbit $(G \times G) \cdot h$ of an element $h \in \P(\End(M))$, where $G \times G$ acts on $\End(M) \cong M^* \otimes M$ via $(g_1,g_2)\cdot(m_1 \otimes m_2)=g_1m_1 \otimes g_2m_2$. Thus, $G$ is embedded in $\P(\End(M))$. Define the wonderful compactification $X$ of $G$ as the closure of $(G \times G) \cdot h$ in $\P(\End(M))$. In characteristic zero, $M$ can be taken as the highest weight module $V(\lambda)$ for any regular weight $\lambda$. The wonderful compactification will be independent of the chosen weight $\lambda$. In positive characteristic, the Steinberg module is such a suitable module.

\begin{prop}[Theorem 3.1 in \cite{CP83}] 
\label{wundervolleEinbettung}
Let $X$ be the wonderful compactification of an adjoint group $G$ of rank $l$.
\begin{aufz}
\item $X$ is smooth.
\item $X \setminus (G \times G) \cdot h$ is the union of $l$ smooth divisors $S_1,\ldots,S_l$, which intersect transversally.
\item There is a bijection between the subsets of $D\dpp\{1,2,\ldots,l\}$ and the set of $G \times G$-orbits in $X$ given by $I \mapsto X_I^\circ$ where the closure of $X_I^\circ$ is $X_I:=\bigcap_{i \notin I} S_i$.
\item The unique closed $G \times G$-orbit $X_\emptyset=\bigcap_{i=1}^l S_l$ is isomorphic to $G/B \times G/B$ and also denoted by $Y$.
\end{aufz}
\end{prop}

So for two subsets $I,J \subset D$ one has $I \subset J \Leftrightarrow X_I \subset X_J$. Every $G \times G$-orbits contains a \textbf{base point} $h_I$, which has the following properties:
\begin{aufz}
\item $(B \times B^-)\cdot h_I$ is dense in $X_I$ and
\item there is a cocharacter $\gamma$ of $T$ such that $h_I=\lim_{t \rightarrow 0} \gamma(t)$.
\end{aufz}

Springer also describes the $B \times B$-orbits and their closure relations explicitly. 

\begin{prop}[Lemma 1.3 in \cite{Sp02}]
\label{BxB-Orbiten}
\begin{aufz}
\item The $B \times B$-orbits in $X$ are
\[
[I,x,w] \dpp (B \times B)\cdot(x,w)\cdot h_I,
\]
where $I \subseteq D$, $x \in W^I$, and $w \in W$. Here $W$ is the Weyl group of $G$, $W_I$ the parabolic subgroup of $W$ generated by the simple reflections corresponding to the roots $\alpha_i$ such that $i \in I$, and $W^I$ is the set of minimal representatives in $W/W_I$. 
\item The relation between the set $\{(I,x,w) \mid I \subseteq D, \ x \in W^I, \ w \in W\}$ and the set of $B \times B$-orbits in $X$ is a bijection.
\item The dimension of a $B \times B$-orbit is given by
\[
\dim\ [I,x,w] = l(w_0) - l(x) + l(w) + \betr{I}.
\]
Here $w_0$ is the longest element in $W$ and $l(w)$ denotes the length of the element $w$ of the Weyl group.
\end{aufz}
\end{prop}

There is a partial order on the set of $B \times B$-orbits defined by 
\[
[I,x',w'] \leq [J,x,w] \ :\Leftrightarrow \ \overline{[I,x',w']} \subseteq \overline{[J,x,w]}.
\]

The closures of the $B \times B$-orbits $[\emptyset,x,w]$, where $x,w \in W$, are the well known Schubert varieties in $Y \cong G/B \times G/B$. Denoting the Schubert variety in $G/B$ corresponding to the Weyl group element $w$ by $S(w)$, one has
\[
\overline{[\emptyset,x,w]} \cong S(xw_0) \times S(w).
\]
The closures of the $B \times B$-orbits $[D,\id,w]$ where $w \in W$ are called \textbf{large Schubert varieties} and denoted by $X(w)$.

\begin{prop}[Proposition 2.4 in \cite{Sp02}]
\label{Abschlussrelation}
Let $[I,x',w']$, $[J,x,w]$ be two $B \times B$-orbits in $X$ where $I,J \subset D$, $x' \in W^I$, $x \in W^J$, and $w',w \in W$. The relation $[I,x',w'] \leq [J,x,w]$ holds if and only if
\begin{aufz}
\item $I \subseteq J$ and
\item there exist $u \in W_I$ and $v \in W_J \cap W^I$ such that
\begin{aufz}
\item $l(wv)=l(w)+l(v)$,
\item $x' \geq xvu^{-1}$, \ and
\item $w'u \leq wv$.
\end{aufz}
\end{aufz}
Here $\leq$ denotes the Bruhat order on the Weyl group.
\end{prop}

\begin{prop}[Lemma 2.4 in \cite{St87}]
The restriction map $\Pic(X) \rightarrow \Pic(Y)$ is injective. Its image is the sub-lattice $\{(-w_0\lambda,\lambda) \mid \lambda \in \Lambda\} \cong \Lambda$.
\end{prop}

So we get $\Lambda \cong \Pic(X)$ via $\lambda \mapsto \Lb{\lambda}$, where $\Lb{\lambda}$ is the line bundle on $X$ which restricts to $\Lb{-w_0\lambda}^{(G/B)} \boxtimes \Lb{\lambda}^{(G/B)}$ on $Y \cong G/B \times G/B$. Here $\Lb{\lambda}^{(G/B)}$ denotes the line bundle $G \times_B k_{-\lambda}$ on $G/B$ associated to the weight $\lambda$.

\begin{prop}[Corollary 8.2 in \cite{CP83} and section 1 in \cite{BP00}]
For every $i=1,\ldots,l$ there is a $G \times G$-invariant section $\sigma_i \in H^0(X,\Lb{{\alpha}_i})$ with divisor $S_i$, that is unique up to multiplication by a scalar.
\end{prop}

\begin{notation}
Let $\lambda,\mu \in \Lambda$ be two weights. Write $\mu \leq \lambda$ if there are non-negative numbers $n_1,...,n_l \in \N_0$ such that $\lambda - \mu = \sum_{i=1}^l n_i {\alpha}_i$. Denote such a collection of numbers by
\[
\vec{n}=(n_1,\ldots,n_l) \in \N_0^l.
\]
The norm of the vector $\vec{n}$ is $\betr{\vec{n}}=\sum_{i=1}^l n_i$. Denoting $\vec{\alpha}=({\alpha}_1,\ldots,{\alpha}_l)$, the scalar product is $\vec{n}\vec{\alpha} = \sum_{i=1}^l n_i {\alpha}_i$. There are two convenient ways to refer to products of $\sigma_i$. Write
\[
\sigma^{\vec{n}} = \sigma_1^{n_1}\cdots\sigma_l^{n_l} = \sigma^{(\lambda - \mu)}
\]
in case $\mu \leq \lambda$ and $\lambda - \mu = \vec{n}\vec{\alpha}$. Write in this case also $\betr{\lambda - \mu}=\betr{\vec{n}}$.
\end{notation}

\begin{prop}[Theorem 8.3 in \cite{CP83}]
\label{ungleichNull}
Let $\lambda \in \Lambda$.
\[
H^0(X,\Lb\lambda) \neq 0 \ \Leftrightarrow \ \exists \mu \in \Lambda^+ \mbox{ such that } \mu \leq \lambda
\]
\end{prop}

\section{Standard monomials for $G/B$}

In this section, we briefly recall the facts we will need about the path model and standard monomials. We refer to \cite{Li98} for a more detailed exposition.

One of the basic results of standard monomial theory is a set $\{p_\pi^{(\lambda)} \mid \pi \in B_\lambda\}$ which forms a basis of the $G$-module $H^0(G/B,\Lb\lambda)$ and has some nice properties. Its elements $p_\pi^{(\lambda)}$ are called \textbf{path vectors} and are indexed by the set $B_\lambda$ of LS-paths $\pi$ of the form $\lambda \in \Lambda^+$. Two of their properties will be of importance later on. First, the end point $\pi(1)$ of a path $\pi \in B_\lambda$ is in $\Lambda$, and the corresponding path vector $p_\pi^{(\lambda)}$ is a weight vector of weight $-\pi(1)$. The second property is that there exists a map $i: B_\lambda \rightarrow W$ which assigns to an LS-path $\pi$ an element of the Weyl group $i(\pi)$ called its initial direction. This is used to define for a path the notion of being standard on a Schubert variety in $G/B$.

\begin{defi}
Let $Y=\bigcup X(\tau_i)$ be a union of Schubert varieties $X(\tau_i)$ in $G/B$. The LS-path $\pi \in B_\lambda$ is  \textbf{standard on $Y$} if and only if $i(\pi) \leq \tau_i$ for at least one $\tau_i$. Here $i(\pi) \in W$ is the initial direction of the path $\pi$ and $\leq$ denotes the Bruhat order on the Weyl group. The associated path vector $p_\pi \in H^0(G/B,\Lb\lambda)$ is \textbf{standard on $Y$} if and only if the path $\pi$ is. 
\end{defi}

\begin{prop}[Theorem 5.3, Corollary 5.2, Theorem 8.6 in \cite{Li98}] \ \ 
\begin{enumerate}
\item The set of path vectors $\{p_\pi \mid \pi \in B_\lambda\}$ of the form $\lambda$ is a basis of $H^0(G/B,\Lb\lambda)$.
\item $\{p_\pi\ueber_Y \mid \pi \in B_\lambda, \ \pi \mbox{ standard on } Y\}$ is a basis of $H^0(Y,\Lb\lambda\ueber_Y)$.
\item $\{p_\pi \mid \pi \in B_\lambda, \ \pi \mbox{ not standard on } Y\}$ is a basis of the kernel of the restriction map $H^0(G/B,\Lb\lambda) \rightarrow H^0(Y,\Lb\lambda\ueber_Y)$.
\end{enumerate}
\end{prop}

\section{Standard monomials for $X$}
\label{BasisvonX}

In  \cite{CM03}, Chiriv\`{i} and Maffei construct standard monomials for the wonderful embedding of a symmetric space $G/H$ that extend the classical standard monomials. The construction and properties of the extended monomials in the case $X=\overline{G}$ are given in this section.

Let $X$ be the wonderful compactification of the group $G$ and $\lambda \in \Lambda^+$ a dominant weight. The line bundle $\Lb\lambda \in \Pic(X)$ can be $\tilde{G} \times \tilde{G}$-linearized which gives a linear action on $H^0(X,\Lb\lambda)$. For any closed $G \times G$-stable subvariety $Z$ of $X$, we consider $H^0(Z,\Lb\lambda)$ as a $\tilde{G} \times \tilde{G}$-module. 

\begin{prop}[Corollary 1.8 in \cite{CM03}]
For all $\lambda \in \Lambda^+$ the restriction $H^0(X,\Lb\lambda) \rightarrow H^0(Y,\Lb\lambda)$ is surjective. 
\end{prop}

Let $\lambda \in \Lambda^+$ be a dominant weight. For any dominant $\mu \leq \lambda$ and any LS-path $\pi \in B_\mu$ there exists a path vector $p_\pi^{(\mu)} \in H^0(Y,\Lb\mu)$. Choose for every $p_\pi^{(\mu)}$ an arbitrary continuation $x_\pi^{(\mu)} \in H^0(X,\Lb\mu)$ such that $x_\pi^{(\mu)}\ueber_Y = p_\pi^{(\mu)}$.

\begin{prop}[Theorem 3.3 in \cite{CM03}]
The set $\M^{(\lambda)} \dpp \{ \sigma^{(\lambda - \mu)} x_\pi^{(\mu)} \mid \mu \in \Lambda^+, \ \mu \leq \lambda, \ \pi \in B_\mu \}$ is a basis of $H^0(X,\Lb\lambda)$.
\end{prop}

The next proposition shows that this basis is compatible with the $G \times G$-orbits in $X$.

The closure $X_I$ of the $G \times G$-orbit corresponding to $I \subseteq D$ satisfies $X_I = \cap_{i \notin I} S_i$. Hence, the restriction of the $G \times G$-invariant section $\sigma_i$ to $X_I$ is non-zero if and only if $i \in I$. In particular, on the unique closed orbit $X_\emptyset=Y$ we have $\sigma_1\ueber_Y=\ldots=\sigma_l\ueber_Y=0$.

\begin{defi}
Let $X_I$ be the closure of a $G \times G$-orbit in $X$ and $\lambda,\mu \in \Lambda^+$ two dominant weights such that $\mu \leq \lambda$. Then $\lambda - \mu=\sum_{i=1}^l n_i {\alpha}_i$ where $n_i \in \N_0$ for all $i\in D$. The monomial $\sigma^{(\lambda - \mu)} x_\pi^{(\mu)} \in \M^{(\lambda)}$ is \textbf{standard on $X_I$} if and only if $n_i=0$ for all $i \notin I$.
\end{defi}

\begin{prop}[Corollary 3.4 in \cite{CM03}]
Let $X_I$ be the closure of a $G \times G$-orbit in $X$. The set
\[
\M_{X_I}^{(\lambda)} \dpp \{ \sigma^{(\lambda - \mu)} x_\pi^{(\mu)}\ueber_{X_I} \ \mid  \sigma^{(\lambda - \mu)} x_\pi^{(\mu)} \mbox{ standard on } X_I \}
\]
is a basis of $H^0(X_I,\Lb\lambda)$.
\end{prop}


\section{A basis for $H^0(Z,\Lb\lambda)$}

Let $X$ be the wonderful compactification of the group $G$ and $Z$ the closure of a $B \times B$-orbit in $X$. In this part, for any $\lambda \in \Lambda^+$ we will define a subset $\M_Z^{(\lambda)}$ of $\M^{(\lambda)}$ which is a basis of the $\tilde{B} \times \tilde{B}$-module $H^0(Z,\Lb\lambda)$.

Consider a $B \times B$-orbit $[I,x,w]$, $I \subseteq D$, $x \in W^I$, $w \in W$, and its closure $Z$. As the intersection $Z \cap Y$ with the unique closed $G \times G$-orbit $Y=X_\emptyset$ is closed and $B \times B$-stable, it is a union of Schubert varieties in $Y \cong G/B \times G/B$. A basis of $H^0(Z \cap Y,\Lb\mu\ueber_{Z \cap Y})$ for $\mu \in \Lambda^+$ is given by the set of restrictions of the standard monomials $p_\pi^{(\mu)}\ueber_{Z \cap Y}$, that are standard on $Z \cap Y$.

For a dominant weight $\mu \in \Lambda^+$ and an LS-path $\pi \in B_\mu$ let $x_\pi^{(\mu)} \in H^0(X,\Lb\mu)$ be an arbitrary continuation of the standard monomial $p_\pi^{(\mu)} \in H^0(Y,\Lb\mu)$ to $X$.

\begin{defi}
\label{standardaufZ}
The path $\pi \in B_\mu$ and $x_\pi^{(\mu)}$ are \textbf{standard on $Z$} if and only if $p_\pi^{(\mu)}= x_\pi^{(\mu)}\ueber_Y$ is standard on the union of Schubert varieties $Z \cap Y$.
\end{defi}

\begin{satz}
\label{linunabh}
Let $Z$ be the closure of the $B \times B$-orbit $[I,x,w]$ in $X$, where $I=\{i_1,\ldots,i_r\} \subseteq D$, $x \in W^I$, and $w \in W$. 
\[
\M_Z^{(\lambda)} \dpp
\left\{ \begin{array}{l|l}
\sigma_{i_1}^{n_1}\cdots\sigma_{i_r}^{n_r} x_\pi^{(\mu)}\ueber_Z  &
\begin{array}{l}
\mu=\lambda-\sum_{k=1}^r n_k \alpha_{i_k} \in \Lambda^+,\ 
n_1,\ldots,n_l \in \N_0,\!\!\! \\
\pi \in B_\mu \mbox{ standard on } Z 
\end{array}
\end{array} \right\}
\]
is a linearly independent subset of $H^0(Z,\Lb\lambda \ueber_{Z})$.
\end{satz}

\begin{beweis}
Consider the equation
\begin{equation}\label{Ansatz}
\sum_{\mu \leq \lambda} \sum_{\pi \in B_\mu} \beta_\pi^{(\mu)} \sigma_{i_1}^{n_1}\cdots\sigma_{i_r}^{n_r} x_\pi^{(\mu)}\ueber_Z = 0,
\end{equation}
where $\beta_\pi^{(\mu)} \in k$ and $x_\pi^{(\mu)}$ is standard on $Z$.
Here the sum is over all $\mu=\lambda-\sum_{k=1}^r n_k \alpha_{i_k} \in \Lambda^+$ such that $n_{1},\ldots,n_{r} \in \N_0$.

All $\sigma_i$ are zero on the closed $G \times G$-orbit $Y$, so the restriction of \eqref{Ansatz} to $Z \cap Y$ yields
\[
\sum_{\pi \in B_\lambda} \beta_\pi^{(\lambda)} x_\pi^{(\lambda)}\ueber_{Z \cap Y} = 0.
\]
But $p_\pi^{(\lambda)}=x_\pi^{(\lambda)}\ueber_Y$ are standard and in particular linearly independent on $Z \cap Y$, so $\beta_\pi^{(\lambda)} =0$ for all $\pi \in B_\lambda$.

The next aim is to show that equation \eqref{Ansatz} also implies $\beta_\pi^{(\mu)}=0$ for all $\pi \in B_\mu$ and $\mu \in M\dpp\{\mu \leq \lambda \mid \mu=\lambda-\sum_{\alpha_{i} \in I} n_i \alpha_{i} \in \Lambda^+,\ n_i\in \N_0\}$, $\mu \neq \lambda$. Therefore, a lexicographic order on the set $\{\mu \in \Lambda^+ \mid \mu \leq \lambda\}$ is defined as follows:

Let $\mu, \mu' \leq \lambda$ where $\mu=\lambda-\sum_{i=1}^l n_i \alpha_i \in \Lambda^+$, $\mu'=\lambda-\sum_{i=1}^l n_i' \alpha_i \in \Lambda^+$. We have $\mu >_{\lex} \mu'$ if and only if there exists $j \in \{1,\ldots,l\}$ such that $n_i=n_i'$ for all $i<j$ and $n_j<n_j'$.

This definition yields a total ordering on the set $\{\mu \in \Lambda^+ \mid \mu \leq \lambda\}$ which can be restricted to the subset $M$ that contains those $\mu = \lambda - \sum_{i \in I} n_i \alpha_i \in \Lambda^+$ where $n_i=0$ for $i \notin I$.

Now consider a weight $\nu \in M$, $\nu < \lambda$ and assume $\beta_\pi^{(\mu)}=0$ for all $\pi \in B_\mu$ and $\mu >_{\lex} \nu$. It remains to show that $\beta_\pi^{(\nu)}=0$ for all $\pi \in B_\nu$. Let $\lambda-\nu=\sum_{k=1}^r n_k \alpha_{i_k}$ and $j \in \{1,\ldots,r\}$ such that $n_{1}=\ldots=n_{j-1}=0$ and $n_{j} \neq 0$.

On the closure $X_{\{\alpha_{i_j},\ldots,\alpha_{i_r}\}}$ of the $G \times G$-orbit $X^\circ_{\{\alpha_{i_j},\ldots,\alpha_{i_r}\}}$ we have $\sigma_{i_1}=...=\sigma_{i_{j-1}}=0$. So if we restrict the sections, equation \eqref{Ansatz} becomes
\[
\sum_{\begin{array}{c} \scriptstyle (m_j,\ldots,m_r) \in \N^{r-j+1},\\[-0,5ex] \scriptstyle \mu = \lambda - \sum_{k=j}^r m_{k}\alpha_{i_k} \in \Lambda^+ \end{array}} \sum_{\pi \in B_\mu} \beta_\pi^{(\mu)} \sigma_{i_j}^{m_j}\cdots \sigma_{i_r}^{m_{r}} x_\pi^{(\mu)}\ueber_{Z \cap X_{\{\alpha_{i_j},\ldots,\alpha_{i_r}\}}} = 0
\]
where $m_j \geq n_j$. All $\sigma_{i_j}^{m_j}\cdots \sigma_{i_r}^{m_{r}} x_\pi^{(\mu)}\ueber_{Z \cap X_{\{\alpha_{i_j},\ldots,\alpha_{i_r}\}}}$ lie in the image of
\[
H^0(Z \cap X_{\{\alpha_{i_j},\ldots,\alpha_{i_r}\}}, \Lb{\lambda-n_{j}\alpha_{i_j}})
\xrightarrow{\sigma_{i_j}^{n_j} \cdot} 
H^0(Z \cap X_{\{\alpha_{i_j},\ldots,\alpha_{i_r}\}}, \Lb{\lambda}).
\]

The wonderful compactification of a group is in particular a complete regular $G \times G$-variety, so Theorem 1.4 in \cite{B98} can be applied. Part (ii) implies that the intersection of an irreducible component of $Z \cap X_J$ with the $G \times G$-orbit $X_J^\circ$ for $J \subseteq I$ is non-empty. Since $\sigma_{i_j}$ is $G$--invariant and does not vanish on the $G \times G$-orbit $X_J^\circ$ for $\alpha_{i_j} \in J$, the multiplication by $\sigma_{i_j}^{n_j}$ is an injective map. Therefore, the set
$\{ \sigma_{i_j}^{m_j}\cdots \sigma_{i_r}^{m_{r}} x_\pi^{(\mu)}\ueber_{Z \cap X_{\{\alpha_{i_j},\ldots,\alpha_{i_r}\}}}\}$ 
is linearly independent if and only if the set of pre-images 
\[
\{ \sigma_{i_j}^{m_j-n_j}\sigma_{i_{j+1}}^{m_{j+1}}\cdots \sigma_{i_r}^{m_{r}} x_\pi^{(\mu)}\ueber_{Z \cap X_{\{\alpha_{i_j},\ldots,\alpha_{i_r}\}}}\} 
\subseteq 
H^0(Z \cap X_{\{\alpha_{i_j},\ldots,\alpha_{i_r}\}}, \Lb{\lambda-n_j\alpha_{i_j}})
\]
is, too. Restricting to the closure of the $G \times G$-orbit $X_{\{\alpha_{i_{j+1}},\ldots\alpha_{i_r}\}}$ the equation
\[
\sum_{\begin{array}{c} \scriptstyle (m_j,\ldots,m_r) \in \N^{r-j+1},\\[-0,5ex] \scriptstyle \mu = \lambda - \sum_{k=j}^r m_{k}\alpha_{i_k} \in \Lambda^+ \end{array}} \sum_{\pi \in B_\mu} \beta_\pi^{(\mu)} \sigma_{i_j}^{m_j-n_j}\sigma_{i_{j+1}}^{m_{j+1}}\cdots \sigma_{i_r}^{m_{r}} x_\pi^{(\mu)}\ueber_{Z \cap X_{\{\alpha_{i_j},\ldots,\alpha_{i_r}\}}} = 0
\]
becomes
\[
\sum_{\begin{array}{c} \scriptstyle (m_{j+1},\ldots,m_r) \in \N^{r-j},\\[-0,5ex] \scriptstyle \mu = \lambda - n_j\alpha_{i_j} - \sum_{k=j+1}^r m_{k}\alpha_{i_k} \in \Lambda^+ \end{array}} \sum_{\pi \in B_\mu} \beta_\pi^{(\mu)} \sigma_{i_{j+1}}^{m_{j+1}}\cdots \sigma_{i_r}^{m_{r}} x_\pi^{(\mu)}\ueber_{Z \cap X_{\{\alpha_{i_{j+1}},\ldots\alpha_{i_r}\}}} = 0,
\]
where $m_{j+1} \geq n_{j+1}$ by induction hypothesis. Hence, the last steps can be repeated with $j+1,j+2,\ldots,r$. Finally, equation \eqref{Ansatz} gives
\[
\sum_{\pi \in B_\nu} \beta_\pi^{(\nu)} x_\pi^{(\nu)}\ueber_{Z \cap Y} = 0.
\] 
As all $x_\pi^{(\nu)}$ are standard on $Z \cap Y$, this implies $\beta_\pi^{(\nu)}=0$ for all $\pi \in B_\nu$.
\end{beweis}

\begin{bem}
The same proof also works for corresponding sets when $X=\overline{G/H}$ is the wonderful compactification of an arbitrary symmetric space. In this case, standard monomials on the closure $Z$ of a $B$-orbit can be defined exactly in the same way. But while the set $\M_Z^{(\lambda)}$ is still linearly independent, it is in general not a basis of $H^0(Z,\Lb\lambda)$. 
\end{bem}

To show that $\M_Z^{(\lambda)}$ is a basis, the dimension of the $\tilde{B} \times \tilde{B}$-module $H^0(Z,\Lb\lambda)$ has to be calculated. For this, we follow the approach of Brion and Polo in \cite{BP00} and generalize some of their results on large Schubert varieties to arbitrary $B \times B$-orbit closures.

First, we need some facts that are obtained by using the methods of Frobenius splitting. Let the algebraic group $G$ be defined over an algebraically closed field of positive characteristic. In \cite{BP00}, Brion and Polo construct a Frobenius splitting $\sigma \in \Hom_{\oo_X}(F_*\oo_X,\oo_X)$ that splits $X$ compatibly with the closures of all $G \times G$-orbits $X_I$ and with the large Schubert varieties $X(w)$ (Theorem 2 in \cite{BP00}). He and Thomsen show in \cite{HT05} that $\sigma$ splits $X$ compatibly also with all $B \times B$-orbit closures (Proposition 7 in \cite{HT05}). From this, the following corollary is obtained.

\begin{folg}
\label{ressur+red}
Let $X$ be the wonderful compactification of the adjoint group $G$ over an algebraically closed field of arbitrary characteristic. Let $\lambda \in \Lambda^+$, $X_I$ the closure of a $G \times G$-orbit and $Z$ the closure of a $B \times B$-orbit in $X$. 
\begin{aufz}
\item \label{ressur+red_Teil1}
The restriction maps $\res_Z:H^0(X,\Lb\lambda) \rightarrow H^0(Z,\Lb\lambda)$ and $\res_{Z \cap Y}:H^0(Z,\Lb\lambda) \rightarrow H^0(Z \cap Y,\Lb\lambda)$ are surjective. Here $\Lb\lambda$ also denotes the restriction of the line bundle to $Z$ respectively $Z \cap Y$.\\
Furthermore, $H^i(Z,\Lb\lambda)=0$ for all $i>0$.
\item The scheme theoretic intersection $Z \cap X_I$ is reduced.
\end{aufz}
\end{folg}

\begin{beweis}
First, consider the case $\car(k)>0$. If $\lambda$ is regular, then the line bundle $\Lb\lambda$ is ample by Lemma 1 in \cite{BP00} and the first assertion follows from Proposition 7.2 in \cite{HT05}. To prove the assertion in case $\lambda$ is not regular, we will use Proposition 1.13 (ii) from \cite{R87}.

In the proof of Theorem 2 in \cite{BP00} it is shown, that $\sigma$ splits $X$ compatibly with the $B \times B$-stable divisor $D^+:=\sum_{i=1}^l X(w_0s_i)$ as well as with the $B^- \times B^-$-stable divisor $D^-:=(w_0,w_0)D^+=\sum_{i=1}^l X^-(s_iw_0)$. Here $B^-$ denotes the opposite Borel of $B$ and $X^-(w)=\overline{B^-wB^-}$, where $w \in W$, is a opposite large Schubert variety in $X$. In particular, $\sigma$ is a $(p-1)D^-$-splitting by Theorem 1.4.10 in \cite{BK04}.

The support $D \dpp \supp((p-1)D^-)=\bigcup_{i=1}^l X^-(s_iw_0)$ contains no $B \times B$-orbit. Indeed, if $x \in D$ such that $(B \times B)x \subseteq D$, then $(B^- \times B^-)(B \times B)x \subseteq D$, because $D$ is $B^- \times B^-$-stable. As $B^-B$ is dense in $G$ and $D$ is closed, this implies $(G \times G)x \subseteq D$. But this is not possible, because $D$ contains no $G \times G$-orbit.
Hence, no irreducible component of a closed union of $B \times B$-orbits is contained in $\supp((p-1)D^-)$. So the $(p-1)D^-$-splitting $\sigma$ splits compatibly with all considered subvarieties $Z$. 

Furthermore, in the proof of Theorem 2 in \cite{BP00} it is shown, that the line bundle associated to the divisor $D^-$ is $\Lb{(p-1)\ro}$. As $\ro$ is a regular weight, this line bundle is ample by Lemma 1 in \cite{BP00}. By the same Lemma, the line bundle $\Lb\lambda$ is generated by its global sections and therefore without base points.

So Proposition 1.13 (ii) in \cite{R87} can be applied and yields the first assertion in positive characteristic. Using Proposition 1.6.2 and Corollary 1.6.3 in \cite{BK04}, this implies the same result in characteristic zero.

The second assertion is an easy consequence of the existence of a splitting in positive characteristic and can be found for example in Proposition 1.2.1 in \cite{BK04}. Corollary 1.6.6 in \cite{BK04} extends the result to characteristic zero. 
\end{beweis}

Recall that $Z$ is the closure of the $B \times B$-orbit $[I,x,w]$ where $I=\{ i_1,\ldots,i_r \} \subseteq D$, $x \in W^I$, and $w \in W$. 

\begin{lemma}
\label{irredKomp}
Let $J \subseteq I$ and $X_J$ be the corresponding $G \times G$-orbit closure. The irreducible components of $Z \cap X_J$ are $\overline{[J,xv,wv]}$ where $v \in W_I \cap W^J$ such that $l(wv)=l(w)+l(v)$.
\end{lemma}

\begin{beweis}
Part (ii) of Theorem 1.4 in \cite{B98} implies that every irreducible component meets the $G \times G$-orbit $X_J^\circ$, so $Z \cap X_J$ is the union of all $\overline{[J,\tilde{x},\tilde{w}]}$ such that $[J,\tilde{x},\tilde{w}] \leq [I,x,w]$. 
\begin{eqnarray*}
[J,\tilde{x},\tilde{w}] \leq [I,x,w] &\Leftrightarrow& \exists v \in W_I \cap W^J \mbox{ such that } l(wv)=l(w)+l(v) \mbox{ and}\\
& & \exists u \in W_J \mbox{ such that } \tilde{x} \geq x v u^{-1} \mbox{ and } \tilde{w}u \leq wv\\
&\Leftrightarrow& \exists v \in W_I \cap W^J \mbox{ such that } l(wv)=l(w)+l(v) \mbox{ and}\\
& & [J,\tilde{x},\tilde{w}] \leq [J,xv,wv]
\end{eqnarray*}
Hence, those $\overline{[J,xv,wv]}$ are the irreducible components.
\end{beweis}

\begin{lemma}
\label{dom}
Let $\mu \in \Lambda$ be a weight.
\[
\mu \notin \bigcap_{i \in I} \ \alpha_i^+ \ \Rightarrow \ H^0(Z \cap Y,\Lb\mu\mid_{Z \cap Y}) = 0
\]
Here for a simple root $\alpha$, denote by $\alpha_i^+$ the set $\{ \lambda \in \Lambda \mid \langle \lambda,\check{\alpha} \rangle \geq 0 \}$.  
\end{lemma}

\begin{beweis}
Let $S(w)$ denote the Schubert variety $G/B$ associated to the element $w$ of the Weyl group. Lemma \ref{irredKomp} implies
\[
Z \cap Y \ = \!
\bigcup_{\begin{array}{c} \scriptstyle v \in W_I\\[-0,5ex] \scriptstyle l(wv)=l(w)+l(v) \end{array}} \overline{[\emptyset,xv,wv]} 
\ \cong \!
\bigcup_{\begin{array}{c} \scriptstyle v \in W_I\\[-0,5ex] \scriptstyle l(wv)=l(w)+l(v) \end{array}} S(xvw_0) \times S(wv). 
\]
The restriction to $Y$ of the line bundle $\Lb\mu$ is the line bundle $\Lb{-w_0\mu}^{(G/B)} \boxtimes \Lb\mu^{(G/B)}$ on $G/B \times G/B$. Dabrowski shows in \cite{D92} that
\[
H^0(S(w),\Lb\mu^{(G/B)}) \neq 0 
\ \Leftrightarrow \ 
\mu \in \alpha^+ \mbox{ for all } \alpha \in \Delta \mbox{ such that } w\alpha \in \Phi^-.
\]
Assume $H^0(Z \cap Y,\Lb\mu) \neq 0$. Then there is a $v \in W_I$ such that $l(wv)=l(w)+l(v)$, $H^0(S(xvw_0),\Lb{-w_0\mu}^{(G/B)}) \neq 0$, and $H^0(S(wv),\Lb{\mu}^{(G/B)}) \neq 0$. We have
\begin{eqnarray*}
H^0(S(wv),\Lb{\mu}^{(G/B)}) \neq 0 
& \Leftrightarrow &
\mu \in \alpha^+ \fa \alpha \in \Delta: wv\alpha \in \Phi^-
\end{eqnarray*}
\begin{eqnarray*}
H^0(S(xvw_0),\Lb{-w_0\mu}^{(G/B)}) \neq 0 
& \Leftrightarrow &
-w_0\mu \in \alpha^+ \fa \alpha \in \Delta: xvw_0\alpha \in \Phi^-\\
& \Leftrightarrow &
-w_0\mu \in (-w_0\alpha)^+ \\
& & \fa \alpha \in \Delta: xvw_0(-w_0\alpha) \in \Phi^-\\
& \Leftrightarrow &
\mu \in \alpha^+ \fa \alpha \in \Delta: xv\alpha \in \Phi^+\\
\end{eqnarray*}

Let $i \in I$. If $v\alpha_i \in \Phi^-$, then $wv\alpha_i \in \Phi^-$, because $l(wv)=l(w)+l(v)$.  If $v\alpha_i \in \Phi^+$, then $xv\alpha_i \in \Phi^+$, because $v \in W_I$ and $x \in W^I$. This shows that $H^0(Z \cap Y,\Lb\mu) \neq 0$ implies $\mu \in \alpha_i^+$ for all $i \in I$.
\end{beweis}

Let $\vec{n}=(n_1,\ldots,n_l) \in \N_0^l$, where $n_i=0$ for all $i \notin I$. Every non-empty open subset $U$ of $Z$ meets the $G \times G$-orbit $X_I^\circ$ on which $\sigma^{\vec{n}}$ is non-zero, so the multiplication by
\[
\sigma^{\vec{n}}: H^0(U,\Lb{-\vec{n}\vec{\alpha}}) \rightarrow H^0(U,\oo_Z)
\]
is an injective map for every $U \subseteq Z$ open. 
Define an ideal sheaf $\sigma^{\vec{n}} \Lb{-{\vec{n}}{\vec{\alpha}}}$ of $\oo_Z$ by
\[
(\sigma^{\vec{n}} \Lb{-{\vec{n}}{\vec{\alpha}}})(U) =
\sigma^{\vec{n}}\ueber_U \cdot \Lb{-{\vec{n}}{\vec{\alpha}}} (U)
\subseteq H^0(U,\Lb{{\vec{n}}{\vec{\alpha}}} \otimes_{\oo(Z)} \Lb{-{\vec{n}}{\vec{\alpha}}}) =
H^0(U,\oo_Z)
\]
for every open set $U$ in $Z$.

\begin{lemma}
\label{Idealgarbe}
Let $\I$ be the ideal sheaf of $Z \cap Y$ in $\oo_Z$. 
\begin{aufz} 
\item 
\label{ersterTeil}
$\I$ is generated by $\sigma_{i_1},\ldots,\sigma_{i_r}$, i.e.
\[
\I = \sum_{\substack{\vec{n} \in \N_0^l \\n_i=0 \fa i \notin I}} \sigma^{\vec{n}} \Lb{-{\vec{n}}{\vec{\alpha}}}.
\]
\item $(\sigma_{i_1},\ldots,\sigma_{i_r})$ form a regular sequence in $\oo_Z$.
\item For all $n \in \N$ we have
\[
\I^n/\I^{n+1} \cong \bigoplus_{\substack{\betr{\vec{n}}=n\\n_i=0 \fa i \notin I}} \sigma^{\vec{n}} \Lb{-{\vec{n}}{\vec{\alpha}}}\ueber_{Z \cap Y}.
\]
\end{aufz}
\end{lemma}

\begin{beweis}
\begin{aufz}
\item 
$\sigma_1,\ldots,\sigma_l$ generate the ideal sheaf $\I_Y$ of $Y$ in $\oo_X$ (see e.g. \cite{BP00} before Corollary 4). All $\sigma_i$ are $G \times G$--invariant, so
\[
\sigma_i\ueber_Z =0 \ \Leftrightarrow \ \sigma_i\ueber_{X_I} =0 \ \Leftrightarrow \ i \notin I.
\]
As the scheme theoretical intersection $Z \cap Y$ is reduced by Corollary \ref{ressur+red}, $\I$ is generated by those $\sigma_i$ where $i \in I$.

\item 
To prove this assertion, the proof of  Corollary 4 in \cite{BP00} can be adapted. For $1 < j \leq r$ define
\[
Z_j \dpp Z \cap X_{\{i_j,\ldots,i_r\}} = Z \cap \bigcap_{k=1}^{j-1} S_{i_k}.
\]
Then we have $\oo_{Z_j} = \oo_Z/(\sigma_{i_1},\ldots,\sigma_{i_{j-1}})$, because $S_{i_k}$ is the divisor corresponding to $\sigma_{i_k}$. Corollary \ref{ressur+red} assures that $Z_j$ is reduced. By Lemma \ref{irredKomp}, the irreducible components of $Z_j$ are $Z_{j,v}\dpp\overline{[J,xv,wv]}$ where $J=\{i_j,\ldots,i_r\}$ and $v \in W_I \cap W^J$ such that $l(wv)=l(w)+l(v)$. As none of these irreducible components is completely contained in $S_{i_j}$, the restriction of $\sigma_{i_j}$ to $Z_{j,v}$ does not vanish for any $v$.

Let $f \in \oo_Z(Z)$ such that $\sigma_{i_j} \cdot f =0$. Then in particular $\sigma_{i_j} \cdot f \ueber_{Z_{j,v}}=0$ holds for the restriction. As $Z_{j,v}$ is irreducible and reduced, $\oo_Z(Z_{j,v})$ is an integral domain. But $\sigma_{i_j}\ueber_{Z_{j,v}} \neq 0$, so $f\ueber_{Z_{j,v}}=0$. This implies $f=0$, and $\sigma_{i_j}$ is no zero divisor in $\oo_{Z_j} = \oo_Z/(\sigma_{i_1},\ldots,\sigma_{i_{j-1}})$.

\item 
Because of \ref{ersterTeil} we get
\[
\I^n = \sum_{\vec{n}} \sigma^{\vec{n}} \Lb{-\vec{n}\vec{\alpha}},
\]
where the sum is over all $\vec{n}=(n_1,\ldots,n_l)\in \N_0^l$ such that $\sum_{i=1}^l n_i=n$ and $n_i=0$ for all $i \notin I$.
From this follows directly
\[
\I^n/\I^{n+1} \cong \bigoplus_{\vec{n}} \sigma^{\vec{n}} \Lb{-{\vec{n}}{\vec{\alpha}}}\ueber_{Z \cap Y}. 
\] \qed
\end{aufz}
\exendbeweis

Using the ideal sheaf $\I$ from the last Lemma, a filtration of the $\tilde{B}\times \tilde{B}$-module $H^0(Z,\Lb\lambda)$ can be defined. Indeed, the $\tilde{B}\times \tilde{B}$-modules
\[
F_n \dpp H^0(Z,\Lb\lambda \otimes \I^n), \mbox{ where } n \in \N_0,
\]
form a finite descending filtration of $H^0(Z,\Lb\lambda)$. For $\vec{n}=(n_1,\ldots,n_l) \in {\N_0}^l$, where $n_i=0$ for all $i \notin I$, the multiplication by 
\[
\sigma^{\vec{n}}: H^0(Z,\Lb{\lambda-\vec{n}\vec{\alpha}}) \rightarrow H^0(Z,\Lb\lambda)
\]
is injective. As $\I$ is generated by those $\sigma_i$ where $i \in I$, and all $\sigma_i$ are invariant under $\tilde{G} \times \tilde{G}$ and therefore in particular under $\tilde{B}\times\tilde{B}$, 
\[
F_{\vec{n}}\dpp\Im(\sigma^{\vec{n}})
\]
is a $\tilde{B} \times \tilde{B}$-submodule of $F_n$, where $n=\betr{\vec{n}}$.

\begin{satz}
\label{Thm7}
\[
F_n=\sum_{\betr{\vec{n}}=n} F_{\vec{n}}
\]
\[
\gr_n H^0(Z,\Lb\lambda) = F_n / F_{n+1} \cong \bigoplus_{\substack{\mu \leq \lambda \mbox{ \scriptsize dom.}\\ \betr{\lambda-\mu}=n}}  H^0(Z \cap Y, \Lb\mu)
\]
\end{satz}

\begin{beweis}
The short exact sequence of sheaves on $Z$ 
\[
0 \rightarrow \Lb\lambda \otimes \I^{n+1} \rightarrow \Lb\lambda \otimes \I^{n} \rightarrow \Lb\lambda \otimes \I^n/\I^{n+1} \rightarrow 0
\]
induces a long exact cohomology sequence
\[
0 \rightarrow F_{n+1} \rightarrow F_n \rightarrow H^0(Z,\Lb\lambda \otimes \I^n/\I^{n+1}) \rightarrow \ldots \, ,
\]
which implies the inclusion
\[
\gr_n H^0(Z,\Lb\lambda) = F_n/F_{n+1} \hookrightarrow H^0(Z,\Lb\lambda \otimes \I^n/\I^{n+1}).
\]
Using the last lemma, we get
\[
H^0(Z,\Lb\lambda \otimes \I^n/\I^{n+1}) =
\bigoplus_{\betr{\vec{n}}=n} \sigma^{\vec{n}} H^0(Z \cap Y,\Lb{\lambda-\vec{n}\vec{\alpha}}).
\]
Here the sum is over all $\vec{n}=(n_1\ldots,n_l)\in \N_0^l$ such that $n_i=0$ for all $i \notin I$. 
Denote $\lambda-\vec{n}\vec{\alpha}$ by $\mu$, and let $j \notin I$. Then 
\[
\langle \mu, \check{\alpha}_j \rangle =
\langle \lambda - \sum_{i \in I} n_i \alpha_i, \check{\alpha}_j \rangle =
\underbrace{\langle \lambda, \check{\alpha}_j \rangle}_{\geq 0 \mbox{\scriptsize, because } \lambda \in \Lambda^+} \!\!
- \, \sum_{i \in I} n_i \underbrace{\langle \alpha_i, \check{\alpha}_j \rangle}_{\leq 0 \mbox{\scriptsize, because } i \neq j}
\geq 0.
\]
So either $\mu \in \Lambda^+ = \bigcap_{i=1}^l \alpha_i^+$ holds or $\mu \notin \bigcap_{i \in I} \alpha_i^+$. But in the second case $H^0(Z \cap Y, \Lb\mu)=0$ by Lemma \ref{dom}. Thus, we get
\[
H^0(Z,\Lb\lambda \otimes \I^n/\I^{n+1}) = \bigoplus_{\substack{\mu \leq \lambda \mbox{ \scriptsize dom.}\\ \betr{\lambda-\mu}=n}} \sigma^{(\lambda-\mu)} H^0(Z \cap Y,\Lb\mu).
\]
Altogether, there is an inclusion
\[
\gr_n H^0(Z,\Lb\lambda) \hookrightarrow \bigoplus_{\substack{\mu \leq \lambda \mbox{ \scriptsize dom.}\\ \betr{\lambda-\mu}=n}} H^0(Z \cap Y,\Lb\mu).
\]

Consider $\vec{n}=(n_1\ldots,n_l)\in \N_0^l$ with $n_i=0$ for all $i \notin I$. As the multiplication by $\sigma^{\vec{n}}: H^0(Z,\Lb{\lambda-\vec{n}\vec{\alpha}}) \rightarrow  H^0(Z,\Lb\lambda)$ is an injective map, its image $F_{\vec{n}}$ is isomorphic to $H^0(Z,\Lb\mu)$ where $\mu=\lambda-\vec{n}\vec{\alpha}$. 
Identifying $\sum_{\betr{\vec{n}}=n} H^0(Z,\Lb{\lambda-\vec{n}\vec{\alpha}})$ with $\sum_{\betr{\vec{n}}=n} F_{\vec{n}} \subseteq H^0(Z,\Lb\lambda)$, we get the well defined restriction map 
\[
\sum_{\betr{\vec{n}}=n} H^0(Z,\Lb{\lambda-\vec{n}\vec{\alpha}}) \rightarrow \bigoplus_{\betr{\vec{n}}=n} H^0(Z \cap Y,\Lb{\lambda-\vec{n}\vec{\alpha}}),
\] 
because any element, that is contained in $\Im(\sigma^{\vec{n}})$ for at least two different $\vec{n}$ such that $\betr{\vec{n}}=n$, restricts to zero on $Y$.
For any dominant weight $\mu \in \Lambda^+$ the restriction map $H^0(Z,\Lb\mu) \rightarrow H^0(Z \cap Y,\Lb\mu)$ is surjective by Corollary \ref{ressur+red}. If $\mu=\lambda-\sum_{i\in I} n_i \alpha_i$ is not dominant, then there is an index $i \in I$ such that $\mu \notin \alpha_i^+$. In this case Lemma \ref{dom} implies $H^0(Z \cap Y,\Lb\mu)=0$. 
This yields the commutative diagram
\[
\xymatrix{\displaystyle 
\displaystyle\sum_{\betr{\vec{n}}=n} F_{\vec{n}} \ar[rr]^\simeq \ar@{^{(}->}[d] &  & \displaystyle\sum_{\betr{\vec{n}}=n} H^0(Z,\Lb{\lambda-\vec{n}\vec{\alpha}})  \ar@{>>}[d] \\
F_n  \ar@{>>}[r] & F_n/F_{n+1} = \gr_n H^0(Z,\Lb\lambda)  \ar@{^{(}->}[r] & \displaystyle\bigoplus_{\substack{\mu \leq \lambda \mbox{ \scriptsize dom.}\\ \betr{\lambda-\mu}=n}} H^0(Z \cap Y,\Lb\mu).\\
}
\]
As this map is surjective, the first part of the assertion follows. 
Now $F_n/F_{n+1} \cong \bigoplus_\mu H^0(Z \cap Y, \Lb\mu)$ implies
\[
F_n \cong \bigoplus_{\mu} H^0(Z \cap Y, \Lb\mu) \oplus F_{n+1}.
\]
As $\sum F_{\vec{n}}$ and $F_{n+1}$ are both submodules of $F_n$ and the above map $\sum F_{\vec{n}} \rightarrow \bigoplus H^0(Z \cap Y, \Lb\mu)$ is surjective, we get
\[
F_n = \sum_{\betr{\vec{n}}=n} F_{\vec{n}} + F_{n+1}.
\]
Iterating the last steps gives 
\begin{eqnarray*}
F_n 
&=& \sum_{\betr{\vec{n}}=n} F_{\vec{n}} + \sum_{\betr{\vec{n}}=n+1} F_{\vec{n}} + F_{n+2}\\
&=& \sum_{\betr{\vec{n}}=n} F_{\vec{n}} + F_{n+2},\\
\end{eqnarray*}
because $F_{\vec{m}} \subseteq F_{\vec{n}}$ if $m=(m_1,\ldots,m_l)$ and $n=(n_1,\ldots,n_l)$ such that $m_i \geq n_i$ for all $1 \leq i \leq l$.
As the filtration is finite, the second part of the assertion follows by induction.
\end{beweis}

\begin{folg}
\label{basis}
The set $\M_Z^{(\lambda)}$ is a basis of $H^0(Z,\Lb\lambda\ueber_Z)$.
\end{folg}

\begin{beweis}
By Proposition \ref{linunabh} the set $\M_Z^{(\lambda)}$ is linearly independent. Theorem \ref{Thm7} implies
\begin{align*}
\dim\ H^0(Z,\Lb\lambda\ueber_Z) 
&= \dim\ \gr \ H^0(Z,\Lb\lambda\ueber_Z) \rb\\
&= \dim\ \bigoplus_{\mu \leq \lambda \mbox{ \scriptsize dom.}} H^0(Z \cap Y, \Lb\mu\ueber_{Z \cap Y}) \rb\\
&= \sum_{\mu \leq \lambda \mbox{ \scriptsize dom.}} \betr{\{ p_\pi^{(\mu)} \mbox{ standard on } Z \cap Y  \}}\\
&= \betr {\M_Z^{(\lambda)}}.  
\end{align*} \qed
\exendbeweis

\section{Standard monomials}

Let $\lambda \in \Lambda^+$ and consider the $\tilde{B} \times \tilde{B}$-module $H^0(X,\Lb\lambda)$. The aim of this section is to construct a basis of this module that -- like classical standard monomials -- has the following properties: 
\begin{aufz}
\item \label{Bed1}
The elements of the basis are indexed by the set of LS-paths $\bigcup_{\mu \leq \lambda \mbox{ \scriptsize dom.}} B_\mu$. We call them path vectors. They are weight vectors whose weight is determined by the end point of the corresponding path. 
\item  \label{Bed2}
Let $Z$ be the closure of a $B \times B$-orbit in $X$. The restriction to $Z$ of those path vectors which are standard on $Z$ with respect to $\lambda$ form a basis of $H^0(Z,\Lb\lambda)$.
\item  \label{Bed3}
Let $Z$ be the closure of a $B \times B$-orbit in $X$. The restriction to $Z$ of those path vectors which are not standard on $Z$ with respect to $\lambda$ form a basis of the kernel of the restriction map $H^0(X,\Lb\lambda) \rightarrow H^0(Z,\Lb\lambda)$.
\end{aufz}

\begin{defi}
\label{bzgl}
Let $\lambda, \mu \in \Lambda^+$ be dominant weights and $Z$ the closure of the $B \times B$-orbit $[I,x,w]$ in $X$. The LS-path $\pi \in B_\mu$ and the corresponding path vector are called \textbf{standard on $Z$ with respect to $\lambda$} if $\pi$ is standard on $Z$ and $\mu \leq \lambda$ such that $\mu=\lambda-\sum_{i=1}^l n_i \alpha_i \in \Lambda^+$ where $n_i=0$ for all $i \notin I$.
\end{defi}

The set $\M^{(\lambda)}=\{\sigma^{(\lambda-\mu)} x_\pi^{(\mu)} \mid \mu \in \Lambda^+,\ \mu \leq \lambda, \ \pi \in B_\mu \}$ defined in section \ref{BasisvonX} is a basis of $H^0(X,\Lb\lambda)$ that is compatible with the restriction to $B \times B$--orbit closures. It fulfils the first two properties. This is true for arbitrary continuations $x_\pi^{(\mu)} \in H^0(X,\Lb\mu)$ of the standard monomials $p_\pi^{(\mu)} \in H^0(Y,\Lb\mu)$.

Now choose for any $\lambda \in \Lambda^+$ and $\pi \in B_\lambda$ an extension $x_\pi^{(\lambda)}$ of the standard monomial $p_\pi^{(\lambda)} \in H^0(Y,\Lb\lambda)$ to $X$. In general, the chosen set $\M^{(\lambda)}=\{\sigma^{(\lambda-\mu)} x_\pi^{(\mu)} \mid \mu \in \Lambda^+,\ \mu \leq \lambda, \ \pi \in B_\mu \}$ does not have property \ref{Bed3}, because the restriction of $\sigma^{(\lambda-\mu)} x_\pi^{(\mu)}$ to a $B \times B$-orbit on which $\pi$ is not standard does not need to be zero. But starting from those, new standard monomials $\sigma^{(\lambda-\mu)} y_\pi^{(\mu)}$ can be constructed as linear combinations of the $\sigma^{(\lambda-\mu)} x_\pi^{(\mu)}$, which have all three properties.

\begin{satz}
\label{Konstruktion der y}
Let $\lambda \in \Lambda^+$ be a dominant weight. For each $\pi \in B_\lambda$ there is a global section $y_\pi^{(\lambda)} \in H^0(X,\Lb\lambda)$ such that $y_\pi^{(\lambda)}\ueber_Y=p_\pi^{(\lambda)}$ and $y_\pi^{(\lambda)}\ueber_Z=0$ for all $B \times B$-orbit closures $Z$ on which $\pi$ is not standard.  
\end{satz}

\begin{beweis}
The claim is proven by constructing $y_\pi^{(\lambda)}$ recursively for all $\lambda$. Let $\lambda \in \Lambda^+$. Assume that $y_\nu^{(\mu)}$, where $\nu \in B_\mu$, is already constructed for all dominant $\mu<\lambda$. Take a path $\pi \in B_\lambda$ and consider
\[
\hat{Z}_\pi \dpp \bigcup_{\substack{\pi \mbox{ \scriptsize not standard}\\ \mbox{\scriptsize on } [I,x,w]}} [I,x,w] \subseteq X. 
\]
$\hat{Z}_\pi$ is closed and $B \times B$-stable. Its irreducible components $Z_1,\ldots,Z_t$ are closures of $B \times B$-orbits. The restriction of $x_\pi^{(\lambda)}$ to each $Z_i$ is a linear combination of elements in $\M^{(\lambda)}_{Z_i}$ of the same weight. Hence, there are coefficients $\alpha_{i\nu}, \beta_{i\nu} \in k$ such that
\[
x_\pi^{(\lambda)}\ueber_{Z_i} = 
\sum_{\substack{\nu \in B_\lambda\\ \nu \mbox{ \scriptsize standard on }Z_i\\ \nu(1)=\pi(1)}} \alpha_{i\nu} x_\nu^{(\lambda)}\ueber_{Z_i} \ + 
\sum_{\substack{\mu<\lambda \mbox{ \scriptsize dom.}\\ \nu \in B_\mu\\ \nu \mbox{ \scriptsize standard on }Z_i \mbox{ \scriptsize w.r.t. } \lambda\\ \nu(1)=\pi(1)}} \beta_{i\nu} \sigma^{(\lambda-\mu)} y_\nu^{(\mu)}\ueber_{Z_i}.
\]
This yields for the restriction to $Z_i \cap Y$ 
\[
\begin{array}{ccc}
x_\pi^{(\lambda)}\ueber_{Z_i \cap Y} &=& \sum \alpha_{i\nu} x_\nu^{(\lambda)}\ueber_{Z_i \cap Y}\\
\parallel & & \parallel \\
p_\pi^{(\lambda)}\ueber_{Z_i \cap Y} &=& \sum \alpha_{i\nu} p_\nu^{(\lambda)}\ueber_{Z_i \cap Y}\\
\end{array}
\]
As $\pi$ is not standard on $Z_i \cap Y$, we have $p_\pi^{(\lambda)}\ueber_{Z_i \cap Y}=0$. But the restrictions $p_\nu^{(\lambda)}\ueber_{Z_i \cap Y}$ form a basis of $H^0(Z_i \cap Y,\Lb\lambda)$, thus $\alpha_{i\nu}=0$ for all $\nu$.

In case $\lambda \in \Lambda^+$ is minimal with respect to the order $\leq$, that means there is no dominant $\mu < \lambda$, this implies that every extension of $p_\pi^{(\lambda)}$ to $X$ has the required properties. Actually, we have an isomorphism $H^0(X,\Lb\lambda) \cong H^0(Z,\Lb\lambda)$, so the choice of the extension is canonical. Denote this extension by $y^{(\lambda)}_\pi$.

For bigger weights $\lambda$, the equation
\[
x_\pi^{(\lambda)}\ueber_{Z_i} = 
\sum_{\substack{\mu<\lambda \mbox{ \scriptsize dom. }\\ \nu \in B_\mu\\ \nu \mbox{ \scriptsize standard on }Z_i \mbox{ \scriptsize w.r.t. } \lambda\\ \nu(1)=\pi(1)}} \beta_{i\nu} \sigma^{(\lambda-\mu)} y_\nu^{(\mu)}\ueber_{Z_i}  
\]
remains. The following argument shows that the coefficients $\beta_{i\nu}$ may be chosen in such a way that $\beta_{i\nu}=\beta_{j\nu}$ for all $i,j \in \{1,\ldots,t\}$. If $\nu$ is not standard on $Z_i$, then $y_\nu^{(\mu)}\ueber_{Z_i}=0$, and $\beta_{i\nu}$ can be chosen arbitrarily. But each $\nu$ is standard on at least one irreducible component $Z_i$. If $\nu$ is standard on two irreducible components $Z_i$ and $Z_j$, then it is standard on their intersection $Z_i \cap Z_j$ as well. This fact is a generalization of the analogous fact for Schubert varieties and it is proved subsequently in Lemma \ref{standardSchnitt}. We have
\[
x_\pi^{(\lambda)}\ueber_{Z_i \cap Z_j} = 
\!\! \sum_{\substack{\nu \mbox{ \scriptsize standard on } Z_i\\ \mbox{ \scriptsize w.r.t. } \lambda}} \beta_{i\nu} \sigma^{(\lambda-\mu)} y_\nu^{(\mu)}\ueber_{Z_i \cap Z_j} = 
\!\! \sum_{\substack{\nu \mbox{ \scriptsize standard on } Z_j\\ \mbox{ \scriptsize w.r.t. } \lambda}} \beta_{j\nu} \sigma^{(\lambda-\mu)} y_\nu^{(\mu)}\ueber_{Z_i \cap Z_j} 
\]
As $\sigma^{(\lambda-\mu)} y_\nu^{(\mu)}\ueber_{Z_i \cap Z_j}=0$ for all $\nu$ that are not standard on $Z_i \cap Z_j$ with respect to $\lambda$, this implies
\[
\sum_{\substack{\nu \mbox{ \scriptsize standard on } Z_i \cap Z_j\\ \mbox{ \scriptsize w.r.t. } \lambda}} \beta_{i\nu} \sigma^{(\lambda-\mu)} y_\nu^{(\mu)}\ueber_{Z_i \cap Z_j} = 
\!\! \sum_{\substack{\nu \mbox{ \scriptsize standard on } Z_i \cap Z_j\\ \mbox{ \scriptsize w.r.t.} \lambda}} \beta_{j\nu} \sigma^{(\lambda-\mu)} y_\nu^{(\mu)}\ueber_{Z_i \cap Z_j}
\]
But all appearing $\sigma^{(\lambda-\mu)} y_\nu^{(\mu)}\ueber_{Z_i \cap Z_j}$ are linearly independent. Therefore we get $\beta_{i\nu}=\beta_{j\nu}$. 
Defining $\beta_{\nu}\dpp\beta_{1\nu}=\ldots=\beta_{t\nu}$ leads to
\[
x_\pi^{(\lambda)}\ueber_{\hat{Z}_\pi} = \sum \beta_\nu \sigma^{(\lambda-\mu)} y_\nu^{(\mu)}\ueber_{\hat{Z}_\pi}.
\]
Now take
\[
y_\pi^{(\lambda)} \dpp x_\pi^{(\lambda)} - \sum \beta_\nu \sigma^{(\lambda-\mu)} y_\nu^{(\mu)}.
\]
This yields for the restriction to $Y$
\[
y_\pi^{(\lambda)}\ueber_Y = x_\pi^{(\lambda)}\ueber_Y -0 = p_\pi^{(\lambda)}
\]
and for the restriction to $\hat{Z}_\pi$
\[
y_\pi^{(\lambda)}\ueber_{\hat{Z}_\pi} = x_\pi^{(\lambda)}\ueber_{\hat{Z}_\pi} - \sum \beta_\nu \sigma^{(\lambda-\mu)} y_\nu^{(\mu)}\ueber_{\hat{Z}_\pi} = 0.
\] 
\end{beweis}

\begin{folg}
Let $\lambda \in \Lambda^+$. The set
\[
\S^{(\lambda)}\dpp\{\sigma^{(\lambda-\mu)} y_\nu^{(\mu)} \mid \mu \in \Lambda^+,\ \mu \leq \lambda,\ \pi \in B_\mu \}
\]
has properties \ref{Bed1}, \ref{Bed2} and \ref{Bed3}. 
\end{folg}

To complete the proof of Proposition \ref{Konstruktion der y}, it remains to show

\begin{lemma}
\label{standardSchnitt}
Let $Z_1$, $Z_2$ be irreducible components of $\hat{Z}_\pi$ where $\pi \in B_\lambda$. If $\mu < \lambda$ is dominant and $\nu \in B_\mu$ is standard on $Z_1$ and $Z_2$ with respect to $\lambda$, then $\nu$ is also standard on $Z_1 \cap Z_2$ with respect to $\lambda$.
\end{lemma}

\begin{beweis}
Let ${Z_1}=\overline{[I_1,x_1,w_1]}$, ${Z_2}=\overline{[I_2,x_2,w_2]}$, and $\lambda-\mu=\sum n_k \alpha_k$. From $\nu \in B_\mu$ standard on $Z_i$ follows $n_k=0$ $\forall k \notin I_i$, hence $n_k=0$ $\forall k \notin I_1 \cap I_2$. As we have $Z_1 \cap Z_2 \subseteq X_{I_1 \cap I_2}$, the weight $\mu$ has the property stated in definition \ref{bzgl}.

By definition, the path $\nu \in B_\mu$ is standard on $Z_i$ if and only if $\nu$ is standard on $Z_i \cap Y$. It remains to show that a path which is standard on two Schubert varieties in $G/B \times G/B$ is also standard on their intersection.

Let $\nu \in B_\mu$ be standard on $Y_i=\overline{[\emptyset,x_i,w_i]}$, $i=1,2$. We claim that in that case $\nu$ is also standard on $Y_1 \cap Y_2$. Indeed, $\nu$ is standard on $Y=\overline{[\emptyset,x,w]} \cong S(xw_0) \times S(w)$ if and only if $i(\nu) \leq (xw_0,w)$, where $i(\nu)$ is the initial direction of the path $i(\nu)$. Denoting $i(\nu)=(\tilde{x}w_0,\tilde{w})$, $\nu$ is standard on $Y$ if $\tilde{x} \geq x$ and $\tilde{w} \leq w$. In particular, $\nu$ is standard on $\overline{[\emptyset,\tilde{x},\tilde{w}]}$.
\begin{eqnarray*}
\nu \mbox{ standard on } Y_1 \mbox{ and } Y_2
& \Leftrightarrow & \tilde{x} \geq x_1,\ \tilde{x} \geq x_2,\ \tilde{w} \leq w_1,\ \tilde{w} \leq w_2\\
& \Leftrightarrow & [\emptyset,\tilde{x},\tilde{w}] \leq [\emptyset,x_1,w_1] \mbox{ and } [\emptyset,\tilde{x},\tilde{w}] \leq [\emptyset,x_2,w_2]\\
& \Leftrightarrow & \overline{[\emptyset,\tilde{x},\tilde{w}]} \subseteq Y_1 \cap Y_2\\
& \Rightarrow & \nu \mbox{ standard on } Y_1 \cap Y_2
\end{eqnarray*} 
\end{beweis}


\end{document}